\documentclass{article}

\usepackage{microtype}
\usepackage{graphicx}
\usepackage{subfigure}
\usepackage{booktabs} 

\usepackage{graphicx}
\usepackage{grffile}
\usepackage{multirow} 
\usepackage{caption}
\usepackage{array}
\usepackage{amsmath}
\usepackage{amsthm}
\usepackage{amsfonts}
\usepackage{amssymb}
\usepackage{mathrsfs}
\usepackage{enumerate}
\usepackage{empheq}

\newcommand{\EE}{\mathbb{E}}

\newcommand{\calO}{\mathcal{O}}

\newcommand{\frakX}{\mathfrak{X}}

\newcommand{\RR}{\mathbb{R}}
\newcommand{\NN}{\mathbb{N}}

\makeatletter
\newcommand*\dotp{\mathpalette\dotp@{.5}}
\newcommand*\dotp@[2]{\mathbin{\vcenter{\hbox{\scalebox{#2}{$\m@th#1\bullet$}}}}}
\makeatother

\newcommand{\dist}{\mathrm{dist}}

\DeclareMathOperator*{\argmin}{arg\,min}

\usepackage{enumerate}
\usepackage[shortlabels]{enumitem} 
\usepackage{cleveref}

\usepackage[accepted]{icml2020}

\newtheorem{definition}{Definition}
\newtheorem{proposition}{Proposition}
\newtheorem{lemma}{Lemma}

\newtheorem{thm}{Theorem}
\newtheorem{coro}{Corollary} 
\newtheorem{assum}{Assumption}

\newcommand{\proj}[1]{\mathrm{proj}_{#1}}

\newcommand{\inner}[2]{\left\langle #1, #2 \right\rangle}

\newcommand{\shaocong}[1]{{#1}} 

\icmltitlerunning{Impact of Model Incoherence on Convergence of SGD with Random Reshuffle}

\begin{document}

\twocolumn[ 
\icmltitle{Understanding the Impact of Model Incoherence on Convergence of Incremental SGD with Random Reshuffle}




\begin{icmlauthorlist}
\icmlauthor{Shaocong Ma}{utah}
\icmlauthor{Yi Zhou}{utah}
\end{icmlauthorlist}
\icmlaffiliation{utah}{Department of Electrical and Computer Engineering, University of Utah, Salt Lake City, United States of America}
\icmlcorrespondingauthor{Yi Zhou}{yi.zhou@utah.edu}
 
\vskip 0.3in
]



\printAffiliationsAndNotice{}  

\begin{abstract}
Although SGD with random reshuffle has been widely-used in machine learning applications, there is a limited understanding of how model characteristics affect the convergence of the algorithm. In this work, we introduce model incoherence to characterize the diversity of model characteristics and study its impact on convergence of SGD with random reshuffle \shaocong{under weak strong convexity}. Specifically, {\em minimizer incoherence} measures the discrepancy between the global minimizers of a sample loss and those of the total loss and affects the convergence error of SGD with random reshuffle.  In particular, we show that the variable sequence generated by SGD with random reshuffle converges to a certain global minimizer of the total loss under full minimizer coherence. The other {\em curvature incoherence} measures the quality of condition numbers of the sample losses and determines the convergence rate of SGD.  With model incoherence, our results show that SGD has a faster convergence rate and smaller convergence error under random reshuffle than those under random sampling, and hence provide justifications to the superior practical performance of SGD with random reshuffle. 
\end{abstract}

\section{Introduction}\label{sec: intro}
We study the following finite-sum optimization problem that covers many important machine learning applications.
\begin{align}
	\min_{\theta\in \mathbb{R}^d} f(\theta) := \frac{1}{n}\sum_{i=1}^{n} \ell_i(\theta), \tag{P}
\end{align}
where $\theta\in \mathbb{R}^d$ corresponds to the model parameters, $f: \mathbb{R}^d\to \mathbb{R}$ denotes the total loss and  each $\ell_i$ corresponds to the sample loss of the $i$-th data sample. Such a problem formulation covers a variety of machine learning problems including support vector machine, logistic regression, matrix completion and neural network training, etc. \shaocong{The common approach is minimizing the error of a predictive model over all data samples in a dataset, an i.i.d. assumption on the data typically decomposes the error into a sum of sample errors.}
 
A standard and widely-applied algorithm that solves the problem (P) is the stochastic gradient descent (SGD) algorithm, which
has been well studied in both convex optimization \cite{bottou2018optimization,robbins_1951,Nemirovski_2009,Lan_2012} and nonconvex optimization \cite{bottou2018optimization,Ghadimi_2016,Ghadimi_2016b}. In these works, the SGD adopts a random sampling with replacement scheme (referred to as {\em random sampling}) and its analysis is based on a bounded stochastic variance assumption. 
 Although such an SGD framework yields theoretically optimal convergence rate \cite{Rakhlin2012}, it cannot fully explain the superior practical performance of SGD in modern machine learning applications where the models are typically over-parameterized and SGD usually adopts the incremental sampling with random reshuffle scheme (referred to as {\em random reshuffle}). Hence, it is desired to develop novel SGD frameworks that provide better understanding of the superior practical performance of SGD with random reshuffle.

Toward this goal, some existing works have proposed various novel analysis frameworks that lead to improved convergence rates of SGD with random sampling. In specific, \cite{tseng1998incremental,solodov1998incremental} introduced a strong growth condition that bounds the maximum sample loss gradient norm in terms of the total loss gradient norm. Under such a condition, \cite{schmidt2013fast} established a sublinear convergence rate and a linear convergence rate for SGD with random sampling in convex and strongly convex optimization, respectively. More recently, \cite{vaswani2018fast} proposed a more relaxed weak growth condition and established similar convergence rate results for SGD with random sampling. In another recent work, \cite{ma2017power} considered an interpolation setting where the model overfits all the data points so that all the sample losses share a unique global minimizer, and they showed that SGD with random sampling achieves a linear convergence rate under strong convexity. On the other hand, another line of works studied SGD with random reshuffle and established sublinear convergence rates under strong convexity, e.g.,  \cite{haochen2018random,nagaraj2019sgd,shamir2016without}. However, the analysis in these works are based on traditional assumptions (e.g., Lipshcitzness, boundedness) that do not emphasize model characteristics, and a sufficiently small step size (typically $\mathcal{O}(n^{-1})$) is required to justify the advantage of random reshuffle over random sampling. 
In particular, these technical settings are not practical in modern machine learning training scenarios where the models are typically {\em over-parameterized} and a {\em constant-level} step size is adopted.
Therefore, it is of great importance and interest to develop a novel theoretical framework for SGD with random reshuffle that characterizes the impact of model characteristics on its convergence under a practical constant step size. In specific, we are interested in studying SGD with random reshuffle in the following aspects.
 \begin{itemize}[topsep=0pt,noitemsep, leftmargin=*]
 	\item The convergence results of SGD with random reshuffle studied in the existing works are {\em in-expectation} with regard to the randomness of reshuffle under a sufficiently small step size. Can we prove {stronger} type of convergence of SGD with random reshuffle under over-parameterized models and a constant step size?
 	\vspace{2pt}
 	\item It has been observed in many practical scenarios that SGD with random reshuffle converges faster than SGD with random sampling under a constant step size. Therefore, the framework that we develop for analyzing SGD with random reshuffle is expected to provide theoretical justifications for this phenomenon.
 	\vspace{2pt}
 	\item The existing theoretical frameworks for analyzing SGD either assume the loss is strongly convex or assume the sample losses share a single global minimizer, both of which rule out many practical machine learning problems that are nonconvex and have multiple global minimizers. Therefore, our framework for analyzing SGD with random reshuffle must cover nonconvex scenarios and allow the existence of multiple global minimizers.
 \end{itemize}

\subsection{Our Contributions}
We analyze the convergence of SGD with random reshuffle under a constant step size by exploiting two notions of model incoherence. In specific, we introduce a minimizer incoherence that measures the discrepancy between the global minimizer of a sample loss and that of the total loss. In particular, full minimizer coherence implies that all the sample losses share a global minimum and hence the model is over-parameterized. We also introduce a curvature incoherence that measures the quality of the condition numbers of the sample losses.
Our theoretical results are in two-fold: 
 \begin{itemize}[topsep=0pt,noitemsep, leftmargin=*]
 	\item We first consider the case of full minimizer coherence where all the sample losses share a set of global minimizers. In such a case, we show that the variable sequence generated by SGD with random reshuffle converges to a certain global minimizer under a constant step size, and therefore the algorithm converges deterministically. Then, under full minimizer coherence and restricted strong convexity, we show that SGD with random reshuffle achieves a linear convergence rate, in which the contraction parameter is determined by the curvature incoherence of the sample losses. Moreover, we establish a linear convergence rate for SGD with random sampling in our framework and prove that SGD achieves a faster linear convergence rate under random reshuffle than that under random sampling, which provides justification to the superior performance of SGD with random reshuffle in training over-parameterized models. We further verify these theoretical results via experiments on over-parameterized neural network training.
 	\vspace{2pt}
 	\item Then, we analyze SGD with random reshuffle in the case of minimizer incoherence where the sample losses do not share any global minimizer. Under a constant step size and restricted strong convexity, we show that SGD with random reshuffle converges to a neighborhood of the global minimizer set at a linear convergence rate. In specific, the convergence rate depends on the curvature incoherence of the sample losses and the convergence error is determined by the minimizer incoherence of the sample losses. We show that the convergence rate of SGD with random reshuffle is faster than that of SGD with random sampling, and the convergence error of SGD is smaller under random reshuffle than that under random sampling. 
 We verify our theoretical results via experiments on nonconvex phase retrieval.
 \end{itemize}

Our analysis shows that the convergence rate and convergence error of SGD with random reshuffle are in the form of geometric mean, whereas those of SGD with random sampling are in the form of arithmetic mean. Therefore, random reshuffle leads to a better convergence statistics for SGD than random sampling.

%
%

\subsection{Related Works}

\paragraph{SGD with random sampling:} Various theoretical frameworks have been developed for analyzing SGD with random sampling. In specific, \cite{schmidt2013fast} exploited the strong growth condition to show that SGD with random sampling achieves a sublinear convergence rate in the convex case and achieves a linear convergence rate in the strongly convex case. \cite{ma2017power} introduced an interpolation setting, in which they showed that SGD with random sampling achieves a linear convergence rate in the strongly convex case. \cite{vaswani2018fast} proposed a relaxed weak growth condition and established a linear convergence rate for SGD with random sampling under strong convexity. \cite{bottou2018optimization} studied SGD with random sampling under a second moment condition. In \cite{gower2019}, they introduced an expected smooth condition and established linear convergence of SGD with random sampling to a neighborhood of the global minimum.

\textbf{SGD with random reshuffle:} It has been noticed that incremental SGD can achieve a faster convergence rate compared to SGD with random sampling in \cite{bottou2009curiously}. The first theoretical analysis was given in \cite{gurbuzbalaban2015random}, where incremental SGD  is shown to outperform SGD with random sampling under a diminishing stepsize. Random reshuffle has been shown to further improve the convergence rate of traditional SGD from $\calO(\frac{1}{k})$ to $\calO(\frac{1}{k^2})$ in the strongly convex case. Then, in more recent works \cite{haochen2018random}, \cite{nagaraj2019sgd}, and \cite{ying2018stochastic}, it was shown that SGD with random reshuffle outperforms SGD with random sampling after finite epochs under a sufficiently small constant step size and strong convexity.  
\section{Introduction to Model Incoherence}
In this section, we introduce two notions of model incoherence.
Recall the finite-sum optimization problem
\begin{align}
\min_{x\in\RR^d} f(\theta) := \frac{1}{n}\sum_{i=1}^n \ell_i(\theta). \tag{P}
\end{align} 
We make the following standard assumption on the existence of solution set of the problem (P). 
\begin{assum}[Existence of solution set]\label{assum: exist}
	Each sample loss $\ell_i, i=1,...,n$ has a solution set $\Theta_i^*\subset \mathbb{R}^d$, on which its global minimum $\ell_i^* >-\infty$ is attained. The total loss $f$ has a solution set $\Theta^*\subset \mathbb{R}^d$, on which its global minimum $f^*>-\infty$ is attained.      
\end{assum}

In general, the solution sets $\{\Theta_i^*\}_{i=1}^n$ of the sample losses can be different from the solution set $\Theta^*$  of the total loss.

We also make the following standard assumptions on the sample losses, where we denote $\proj{A}(x)$ as the Euclidean projection of $x$ onto set $A$. 
\begin{assum}\label{assum: f}
	The problem (P) satisfies:
	\begin{enumerate}[topsep=0pt, noitemsep, leftmargin=*]
		\item The sample losses are $L$-smooth, i.e., $\forall i$ and $x,y\in\mathbb{R}^d$,
		\begin{align}
			\ell_{i}(x) \le \ell_i(y) + \inner{x-y}{\nabla \ell_i(y)} + \frac{L}{2}\|x-y\|^2; \nonumber
		\end{align}
		\item Every sample loss $\ell_i$ \shaocong{is $\mu_i$-weakly strong convex on $\Theta_i^* \cup \Theta^*$}, i.e., 
		\begin{align}
			\ell_i(\omega) \ge \ell_i(x) + \inner{\omega - x}{\nabla \ell_i(x)} + \frac{\mu_i}{2}\|x-\omega\|^2 \nonumber
		\end{align} 
		holds for all $x\in \mathbb{R}^d, \omega \in \proj{\Theta_i^*}(x),  \proj{\Theta^*}(x) $.
	\end{enumerate}
\end{assum}
\shaocong{We note that the weak strong convexity is a weaker condition than the usual strong convexity and covers a wide range of non-convex problems including phase retrieval \cite{Zhou2016b,zhang_2018}, neural networks \cite{zhong_2017,Zhou2017}, low-rank matrix factorization \cite{Tu_2016}, blind deconvolution \cite{Li_2018}, etc. Also, the weak strong convexity implies the restricted strong convexity under an additional convexity condition. We refer to \cite{karimi2016linear} for further discussions.}


\subsection{Minimizer Incoherence}
In this subsection, we introduce minimizer incoherence to measure the discrepancy between the sample loss solution sets $\{\Theta_i^*\}_{i=1}^n$  and the total loss solution set $\Theta^*$. In \Cref{subsec: compare_condition}, we provide a discussion that outlines the connections between the minimizer incoherence and other loss conditions that have been studied in the existing literature.

\begin{definition}[Minimizer incoherence]\label{def: incoherence}
	The minimizer incoherence $\epsilon_i>0$ of every sample loss $\ell_i$ is defined as
	\begin{align*}
	\epsilon_i :=\sup_{\theta \in \Theta^*}   \ell_i(\theta)  -  \ell_i^\ast, \quad i=1,...,n.
	\end{align*} 
\end{definition}

To elaborate, the minimizer incoherence corresponds to the gap between the highest sample loss that is achievable on the total loss solution set and the global minimum of the sample loss. Intuitively, it measures the incoherence between the sample loss solution set $\Theta_i^*$ and the total loss solution set $\Theta^*$. In particular, when minimizer incoherence vanishes, the following inclusion properties of the solution sets hold.
\begin{proposition}[Minimizer coherence]\label{prop: vanish_incoherence}
	The definition of minimizer incoherence implies that
	\begin{enumerate}[topsep=0pt, noitemsep, leftmargin=*]
		\item If $\epsilon_i=0$ for some $i$, then $\Theta^* \subset \Theta_i^*$;
		\item If $\epsilon_i=0$ for all $i$, then $\Theta^* = \cap_{i=1}^n \Theta_i^*$.
	\end{enumerate}
\end{proposition}
In particular, the second item corresponds to the case where we have full minimizer coherence, i.e., all the sample losses share the set of global minimizers $\Theta^*$. This is common in deep learning applications where the models are over-parameterized to overfit all the data samples (hence have full minimizer coherence) and have multiple global minimizers. 
Moreover, our minimizer incoherence generalizes the interpolation condition proposed in \cite{ma2017power}, which requires all the sample losses to share a unique global minimizer under strong convexity. 

\begin{figure}[bth]
	\centering
	\includegraphics[width=0.4\linewidth]{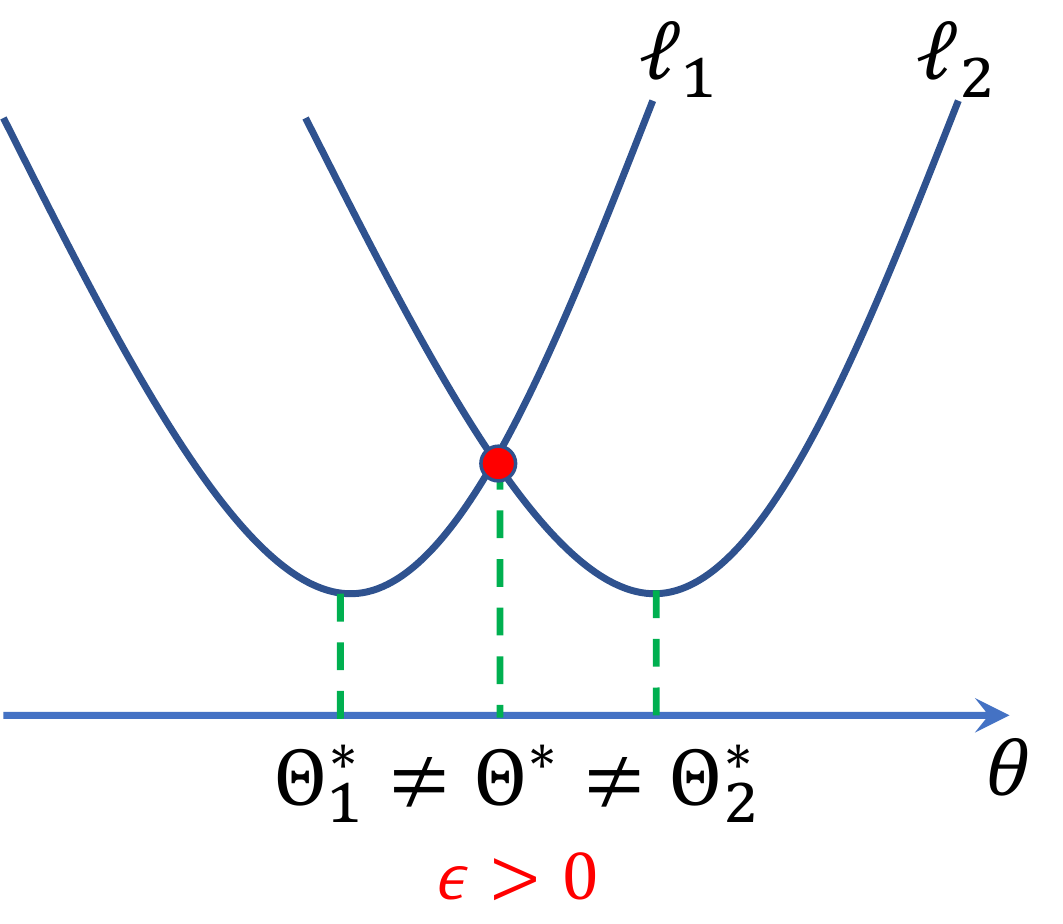}
	\includegraphics[width=0.4\linewidth]{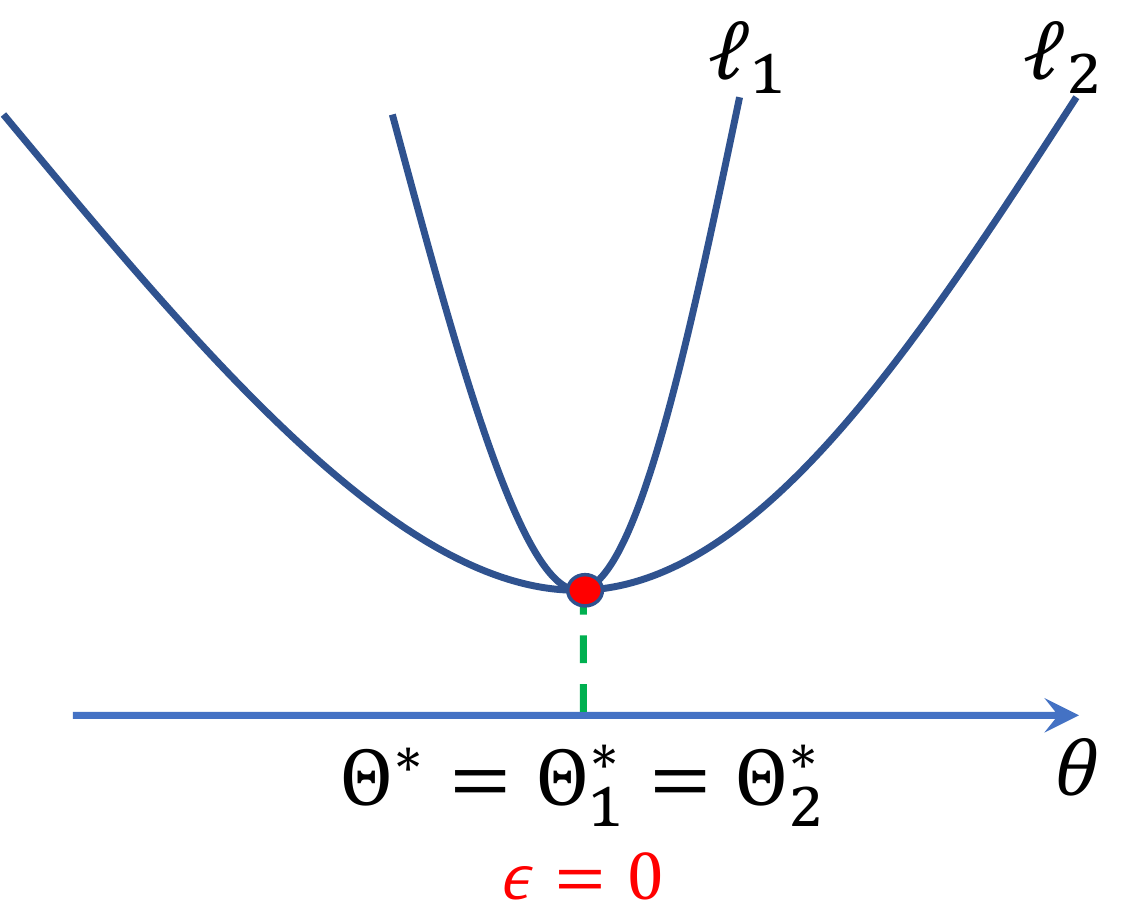}
	\caption{Left: Illustration of minimizer incoherence. Right: Illustration of full minimizer coherence.}\label{fig: 1} 
\end{figure}

\Cref{fig: 1} illustrates the cases of both minimizer incoherence and full minimizer coherence via quadratic sample losses. In fact, many nonconvex machine learning problems have been shown to have either vanishing or small minimizer incoherence, and we provide two illustrative examples below.
\begin{itemize}[topsep=0pt, noitemsep, leftmargin=*]
	\item Phase retrieval \cite{zhang_2018}: In this problem, we take linear measurements of an underlying complex signal $x_0\in \mathbb{C}^d$ with multiple Gaussian vectors $\{a_i \}_{i=1}^n$ and make phaseless observations $y_i = |a_i^\intercal x_0|$. The goal is to recover the complex signal up to a global phase shift by solving the problem
	\begin{align}
		\min_{x\in\mathbb{C}^d} f(x) := \frac{1}{2n}\sum_{i=1}^{n} \big(y_i - |a_i^\intercal x|\big)^2. \nonumber
	\end{align}
	It is clear that all the sample losses share the set of minimizers $\{x_0e^{j\phi}~|~\phi\in (0, 2\pi] \}$ and hence have full minimizer coherence.
	\vspace{2pt}
	\item Over-parameterized neural networks: In deep learning, the neural network model $\theta$ is typically over-parameterized so that the predictor $h_\theta$ can be trained to overfit all the training samples, i.e., $h_\theta(x_i) \approx y_i$ for all $i=1,...,n$. Such overfitting usually achieves a small total loss as well as small sample losses. To justify this, we train a Resnet18 network to overfit the MNIST dataset \shaocong{with the cross-entropy loss. We do not apply any regularization}. \Cref{fig: Resnet18} shows the distribution of the sample losses after $50$ training epochs. One can see that most of the sample losses are below $3\times 10^{-3}$, implying that deep models have very small minimizer incoherence. 

	\begin{figure}[bth]\centering
		\includegraphics[width=0.8\linewidth]{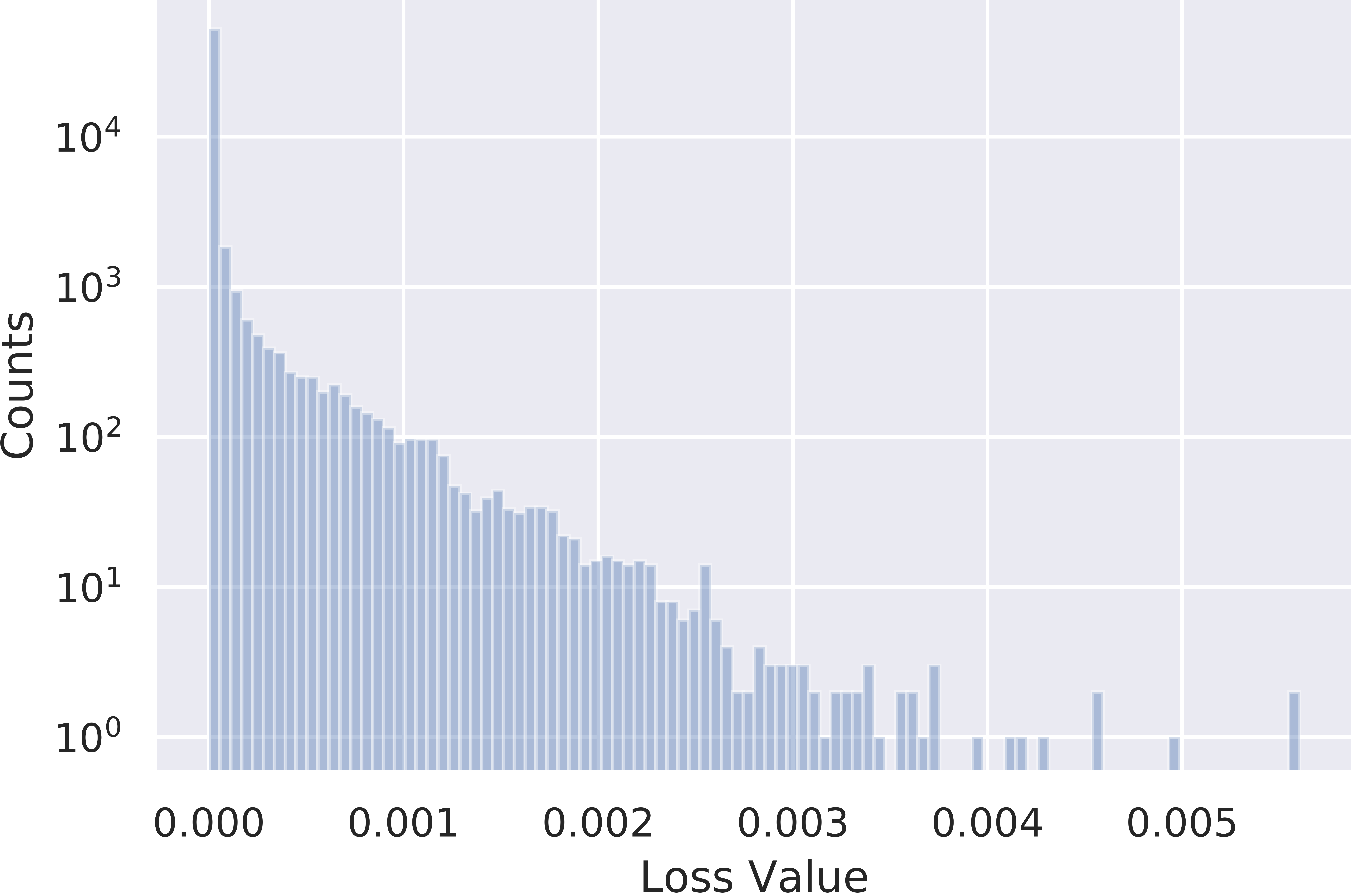} 
		\caption{Distribution of sample losses after training ResNet18 for $50$ epochs on MNIST dataset.} \label{fig: Resnet18}
	\end{figure}
\end{itemize}

\subsection{Curvature Incoherence}


In this subsection, we introduce the curvature incoherence. Recall that the condition number of each sample loss $\ell_i$ is $\frac{L}{\mu_i}$. Then, we define the following curvature incoherence. 


\begin{definition}[Curvature incoherence]\label{def: cur_incoherence}
	The curvature incoherence $\alpha$ of the sample losses is defined as 
	\begin{align*}
	\alpha:= \prod_{i=1}^n\Big(1 - \frac{\mu_i}{L} \Big).
	\end{align*}
In particular, $\alpha$ belongs to the range $[0,1)$.
\end{definition}
As an intuitive understanding, if all the sample losses have a good condition number (i.e., $\frac{L}{\mu_i}\to 1$), then $\alpha$ vanishes and the curvatures of all sample losses are highly coherent.

\section{SGD with Random Reshuffle}
In this section, we introduce the SGD with random reshuffle algorithm and provide some preliminary results on it. 

The SGD algorithm starts with an initialization $\theta_0\in \RR^d$ and applies the following update rule iteratively.
\begin{align}
	\text{(SGD):}~\theta_{k+1} = \theta_k - \eta \nabla\ell_{\xi(k)}(\theta_k), ~k=0,1,...,
\end{align}
where $\eta>0$ is the step size and $\xi(k)$ corresponds to the index of data sample drawn from $\{1,...,n\}$ randomly in the  $k$-th iteration. In this work, we focus on the widely-used incremental sampling with random reshuffle scheme, which is formally defined as follows and is referred to as {\em random reshuffle} for simplicity throughout the paper, . 

\textbf{(Random reshuffle):} {\em In each epoch, we apply a random permutation to the sample indexes, i.e., $\{1,2,...,n\}\overset{\textrm{permute}}{\longrightarrow} \{\xi(0), \xi(1),..., \xi(n-1) \}$. Then, the sample $\xi(k)$ is used in the $k$-th iteration of this epoch.}
%
%

We obtain the following preliminary result for SGD with random reshuffle. 
\begin{lemma}  \label{lemma: 1}
	Let Assumptions \ref{assum: exist} and \ref{assum: f} hold. Apply SGD with random reshuffle  to solve the problem (P) with step size \shaocong{$\eta\leq\frac{1}{L}$}. 
	 Then, the variable sequence $\{\theta_k\}_k$ generated by the algorithm satisfies: for all $k\in \NN$ and all $\omega\in \Theta_{\xi(k)}^*$,
	\begin{align*}
	\|\theta_{k+1} - \omega\|^2 \le \|\theta_{k} -\omega \|^2  - \eta \big(   
	\ell_{\xi(k)}(\theta_{k+1}) -\ell_{\xi(k)}^*
	\big). 
	\end{align*} 
\end{lemma} 
We note that the proof of \Cref{lemma: 1} only requires the sample losses to be restricted convex (i.e., $\mu_i$ can be zero in \Cref{assum: f}.3).
The above lemma characterizes the per-iteration progress of SGD with random reshuffle towards any global minimizer of the sample loss used in the $k$-th iteration. In particular, it implies that $\|\theta_{k+1} - \omega\| \le \|\theta_{k} -\omega \|$, i.e., SGD makes monotonic progress towards the minimizer of the sample loss used in the $k$-th iteration.

\section{Analysis under Full Minimizer Coherence}
In this section, we study the convergence properties of SGD with random reshuffle under full minimizer coherence, i.e., $\epsilon_i=0$ for all $i=1,...,n$ (see \Cref{def: incoherence}) and hence all sample losses share a set of global minimizers $\Theta^*$. 

\subsection{Convergence of SGD Trajectory}
We note that all the results in this subsection only require the sample losses to be restricted convex.
We first characterize the boundedness of the optimization trajectory of SGD with random reshuffle.
\begin{lemma}[Bounded trajectory]\label{lemma: 3}
	Let Assumptions \ref{assum: exist} and \ref{assum: f} hold and assume that the problem (P) has full minimizer coherence. Apply SGD with random reshuffle with step size \shaocong{$\eta\leq \frac{1}{L}$} to solve the problem. Then, the SGD trajectory $\{\theta_k\}_k$ is bounded.
\end{lemma}

As the trajectory of SGD with random reshuffle $\{\theta_k \}_k$ is bounded, it has a compact set of limit points and we denote it as $\frakX$. Also, note that the iteration index sequence  $\{k\}_{k\in \mathbb{N}}$ can be decomposed into $n$ subsequences $\{i(T)\}_T, i=1,...,n$, each of which tracks the SGD iterations that sample the $i$-th data point in the epochs $T=1,2,...$. In particular, we denote $\frakX_i$ as the set of limit points of \shaocong{$\{\theta_{i(T)}\}_T$} and it holds that $\frakX=\bigcup_{i=1}^n \frakX_i$. Moreover, we obtain the following properties regarding the limit point sets of the trajectory of SGD with random reshuffle.

\begin{proposition}[Limit points]\label{prop: 2}
	Under the same conditions as those of \Cref{lemma: 3}, the trajectory of SGD with random reshuffle satisfies the following properties.
	\begin{enumerate}[topsep=0pt, noitemsep, leftmargin=*]
		\item $\frakX_i \subset \Theta_i^*$ for all $i=1,...,n$;
		\item $\frakX_i = \frakX \subset \Theta^*$ for all $i=1,...,n$.
	\end{enumerate}
\end{proposition}

To elaborate, item 1 shows that each sub-trajectory $\{\theta_{i(T)}\}_T$ generated by SGD with random reshuffle is a minimizing sequence for the corresponding sample loss $\ell_i$. Item 2 further strengthens item 1 by showing that  all the sub-trajectories $\{i(T)\}_T, i=1,...,n$ share the same set of limit points, which is a subset of the global minimizer set of the total loss. Intuitively, this is due to the fact that all the sample losses share a set of global minimizers under full minimizer coherence, which guarantees the sub-trajectories of SGD with random reshuffle to have consistent asymptotic properties. 

The proof of \Cref{prop: 2} consists of two major steps. We first exploit full minimizer coherence to prove item 1 and the stationary condition $\|\theta_{k+1} - \theta_k\| \overset{k}{\to} 0$. Then, the stationary condition further guarantees that all sub-trajectories share the same set of limit points and hence implies item 2.

Our main result below further strengthens the convergence properties of the SGD trajectory.
\begin{thm}[Trajectory convergence]\label{thm: conv}
	Under the same conditions as those of \Cref{lemma: 3}, every trajectory $\{\theta_k\}_k$ generated by SGD with random reshuffle converges to a certain global minimizer in $\Theta^*$, i.e., it has a single limit point.
\end{thm}

The above result shows that the entire trajectory of SGD with random reshuffle converges to a certain global minimizer in the case of full minimizer coherence. This implies that full minimizer coherence helps suppress the randomness of the random reshuffle and leads to \shaocong{the point-wise convergence}. Such a deterministic convergence result of SGD with random reshuffle is stronger than other {\em in-expectation} convergence results of SGD with random sampling that are established under various loss conditions (e.g., strong growth condition, interpolation) that imply full minimizer coherence.

\subsection{Convergence Rate Analysis}
%
%
%
%
%
In this subsection, we further study the convergence rate of SGD with random reshuffle under full minimizer coherence. For any point $\theta\in \RR^d$, we denote its distance to an arbitrary set $A\subset \RR^d$ as $\dist_{A}(\theta):=\inf_{u\in A} \|\theta-u\|$.

We obtain the following convergence rate result. 
\begin{thm}[Random reshuffle]\label{thm: rate_coherent-SC}
	Let Assumptions \ref{assum: exist} and \ref{assum: f} hold and assume that the problem (P) has full minimizer coherence. Apply SGD with random reshuffle with step size \shaocong{$\eta = \frac{1}{L}$} to solve the problem. Then, for all epochs $B=1,2,...$, it holds that
	\begin{align}\label{eq: rate1}
		\dist_{\Theta^*}^2(\theta_{nB}) \le \prod_{i=1}^n \Big(1-\frac{\mu_i}{L}\Big)^B \dist_{\Theta^*}^2(\theta_{0}).
	\end{align}
\end{thm}

The above theorem establishes the linear convergence rate of SGD with random reshuffle under full minimizer coherence and the constant step size \shaocong{$\eta = \frac{1}{L}$}. In particular, the convergence rate depends on the curvature incoherence parameter $\alpha = \prod_{i=1}^n (1-\frac{\mu_i}{L})$ (see \Cref{def: cur_incoherence}) that characterizes the quality of the condition numbers of all the sample losses. We also note that in the special case that all the sample losses have the same condition number $\frac{\mu}{L}$, the above convergence rate of SGD with random reshuffle is of order $\calO(1-\frac{\mu}{L})^{nB}$, which meets the convergence rate of full gradient descent under strong convexity.

\subsection{Comparison to Other Sampling Schemes} 
We further analyze the convergence rates of SGD with incremental sampling (i.e., cyclic sampling without random reshuffle) and random sampling under full minimizer coherence and compare them with that of SGD with random reshuffle. 

In fact, under full minimizer coherence, our proof of \Cref{thm: rate_coherent-SC} only rely on the fact that the random reshuffle scheme samples every data point once in each epoch, which is also satisfied by the incremental sampling scheme. Therefore, the convergence rate result in \Cref{thm: rate_coherent-SC} also applies to SGD with incremental sampling and we obtain the following corollary.

\begin{coro}[Incremental sampling]\label{coro: 1}
	Under the same settings as those of \Cref{thm: rate_coherent-SC} and apply SGD with incremental sampling and step size \shaocong{$\eta = \frac{1}{L}$} to solve the problem (P). Then, the convergence rate is also characterized by \cref{eq: rate1}.
\end{coro}

On the other hand, we obtain the following result for SGD with random sampling.


\begin{proposition}[Random sampling]\label{prop: rate_coherent_SGD}
	Under the same settings as those of \Cref{thm: rate_coherent-SC} and apply SGD with random sampling and step size \shaocong{$\eta = \frac{1}{L}$} to solve the problem (P). Then, for all epochs $B=1,2,...$, it holds that
	\begin{align}
	\dist_{\Theta^*}^2(\theta_{nB}) \leq \Big(1 - \frac{\bar{\mu}}{L}\Big)^{nB} \dist_{\Theta^*}^2(\theta_{0}),
	\end{align} 
	where $\bar{\mu} := \frac{1}{n}\sum_{i=1}^n \mu_i$. 
\end{proposition}
The above result establishes a linear convergence rate for SGD with random sampling. Note that the convergence rate depends on the average of the condition numbers of the sample losses. This is different from the convergence rate of SGD with random reshuffle, which depends on the product of the condition numbers of all the sample losses. In particular, by the arithmetic mean-geometric mean (AM-GM) inequality, it holds that
\begin{align}
	\prod_{i = 1}^n \Big( 1 - \frac{\mu_i}{L} \Big) \leq \Big( 1 - \frac{\bar{\mu}}{L} \Big)^n. \label{eq: am-gm}
\end{align}
Therefore, under full minimizer coherence, SGD achieves a faster convergence rate under random reshuffle than that under random sampling. Such a result provides a theoretical justification for the superior performance of SGD with random reshuffle in training over-parameterized models. 

We note that \cite{haochen2018random} also obtains a similar comparison of convergence rate between SGD with random sampling and SGD with random reshuffle. However, their analysis requires the loss to be uniformly strongly convex, whereas our result applies to the broader class of restricted strongly convex functions. Moreover, we established trajectory convergence of SGD with random reshuffle under the existence of multiple global minimizers, whereas their result establishes convergence in expectation under the existence of a unique global minimizer.

\subsection{Empirical Verification}
\begin{figure}[bth]\centering  
	\includegraphics[width=0.49\linewidth]{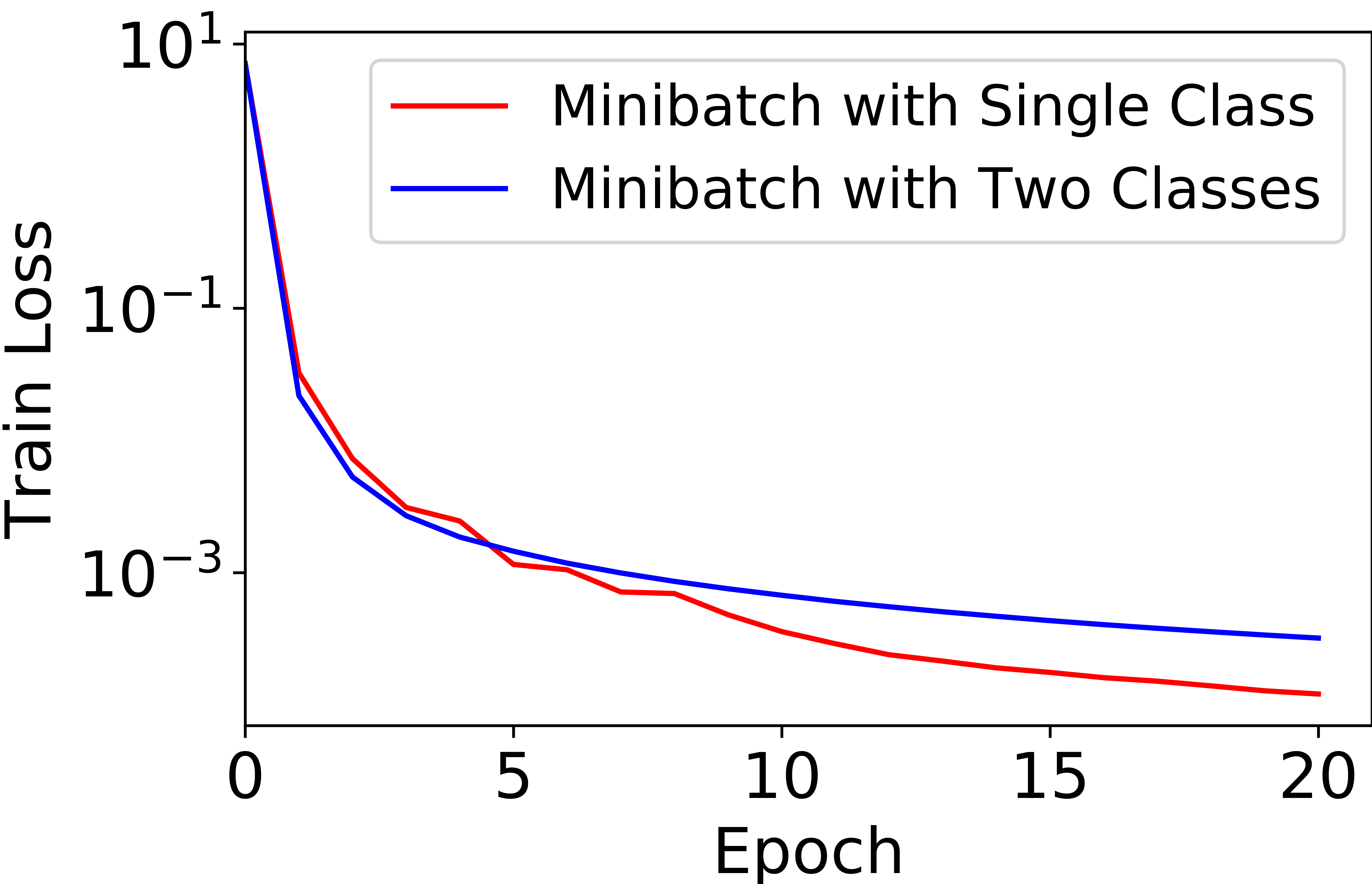} 
	\includegraphics[width=0.49\linewidth]{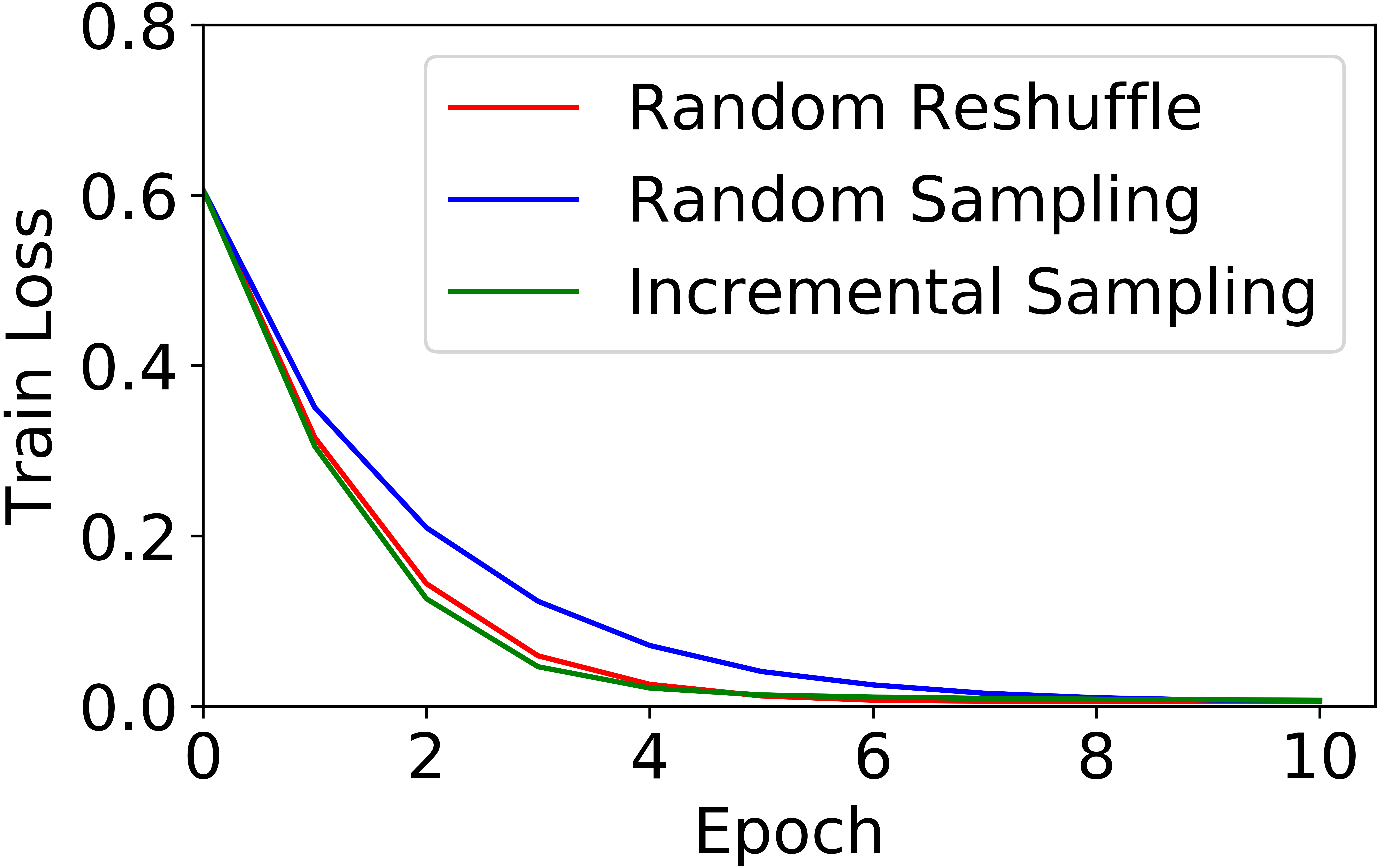}
	\caption{Left: Comparison of convergence of SGD with random reshuffle under different curvature incoherence distributions. Right: Comparison of convergence of SGD under random reshuffle, incremental sampling and random sampling.}\label{fig: 4}
\end{figure}

In this subsection, we verify our theoretical results via experiments. 
We first study the impact of curvature incoherence $\alpha$ on the convergence of SGD with random reshuffle. In specific, we train a Resnet 18 network using SGD with random reshuffle on a mini MNIST dataset that consists of 1000 images of digit ``1'' and 1000 images of digit ``8''. To model different distributions of curvature incoherence, we divide the data samples evenly into 50 fixed mini-batches and consider two different settings: 1) each mini-batch contains 50\% images of digit ``1'' and 50\%  images of digit ``8''; and 2) each mini-batch contains either images of digit ``1'' or images of digit ``8''. In both settings, the average condition numbers of the sample losses are different. 
\Cref{fig: 4} (Left) shows the training loss curves of SGD with random reshuffle in these two settings starting from the same initialization point. It can be observed that SGD with random reshuffle converges faster in the second setting. This implies that the average condition number of the sample losses in the first setting is better than that in the second setting.

Next, we further compare the empirical convergence of SGD under random reshuffle with that under incremental sampling and random sampling. We train a Resnet18 on $4096$ images sampled from CIFAR10 using SGD with the three sampling schemes. We use learn rate $\eta=0.03$, batch-size $128$ and a fixed initialization model that is trained by SGD with random reshuffle for one epoch (with learning rate $0.0025$) using a pre-trained ImageNet model. \Cref{fig: 4} (Right) shows the training loss curves of the three algorithms. It can be seen that SGD with random reshuffle and incremental sampling have a comparable convergence speed, both of which are faster than that of SGD with random sampling. This observation fully supports our theoretical comparison in \cref{eq: am-gm}.

\section{Analysis under Minimizer Incoherence}
In this section, we study the convergence of SGD with random reshuffle under minimizer incoherence where $\epsilon_i > 0$ for some $i\in\{1,...,n\}$. In such a case, the minimizer sets of the sample losses are different from that of the total loss. For simplicity, we assume the minimizer incoherences of all the sample losses are bounded by $\epsilon:= \max_{i\in \{1,...,n\}} \epsilon_i$. 


\subsection{Convergence Rate Analysis}

%
We first show the boundedness of the trajectory of SGD with random reshuffle under minimizer incoherence and a constant step size.
\begin{lemma}[Bounded trajectory]\label{lemma: 2}
	Let Assumptions \ref{assum: exist} and \ref{assum: f} hold and assume that the problem (P) has bounded minimizer incoherence $\epsilon$. Apply SGD with random reshuffle with step size \shaocong{$\eta = \frac{1}{L}$} to solve the problem. Then, the SGD trajectory $\{\theta_k\}_k$ is bounded.
\end{lemma}
The above result generalizes the bounded trajectory result in \Cref{lemma: 3}, which is proved under full minimizer coherence. 
Next, we obtain the following result regarding the convergence rate of SGD with random reshuffle under minimizer incoherence and a constant step size.

\begin{thm}[Random reshuffle]\label{thm: incohe_sc}
	Let Assumptions \ref{assum: exist} and \ref{assum: f} hold. Apply SGD with random reshuffle with step size $\eta=\frac{1}{L}$ to solve the problem (P). Then, for all epochs $B=1,2,...$,
		\begin{align}
		&\mathbb{E} 	\dist_{\Theta^*}^2(\theta_{nB}) \le \alpha^B \Big(	\dist_{\Theta^*}^2(\theta_{0}) \!-\! \frac{2\epsilon \overline{M}}{L(1\!-\!\alpha)} \Big) \!+\! \frac{2\epsilon \overline{M}}{L(1\!-\!\alpha)}, \nonumber 
		\end{align}
{where} $\overline{M} = \mathbb{E}_{\xi}  \bigg[\sum_{k = 0}^{n - 1} \prod_{s=k+1}^{n-1}  \Big(1 - \frac{\mu_{\xi(s)}}{L} \Big) \bigg]$ and $\alpha$ corresponds to the curvature incoherence.
\end{thm}

The above result shows that SGD with random reshuffle converges linearly to a neighborhood of the global minimizer set under minimizer incoherence. Similar to the full minimizer coherence case, the convergence rate coefficient is determined by the curvature incoherence $\alpha$. Moreover, the size of the neighborhood is characterized by the minimizer incoherence $\epsilon$ and the condition numbers of the sample losses. This explains why over-parametrized models such as neural networks can be trained to achieve a small loss by SGD with constant step size: they have very small minimizer incoherence, as demonstrated by the experiment in \Cref{fig: Resnet18}. 
In general, a higher minimizer incoherence and worse condition numbers lead to a larger convergence error of SGD.

\subsection{Comparison to Other Sampling Schemes}
We also obtain the convergence rates of SGD with incremental sampling and random sampling under minimizer incoherence and a constant step size. 

In specific, for SGD with incremental sampling, we denote $\{\sigma(0), \sigma(1),...,\sigma(n-1)\}$ as a specific permutation of the data sample indexes. Such a permutation is fixed throughout the entire training process under incremental sampling. We obtain the following result on SGD with incremental sampling.

\begin{coro}[Incremental sampling]\label{coro: 2}
	Let Assumptions \ref{assum: exist} and \ref{assum: f} hold. Apply SGD with incremental sampling with step size $\eta=\frac{1}{L}$ to solve the problem (P). Then, for all epochs $B=1,2,...$,
	\begin{align}
	&\dist_{\Theta^*}^2(\theta_{nB}) \le \alpha^B \Big(	\dist_{\Theta^*}^2(\theta_{0}) - \frac{2\epsilon \widetilde{M}}{L(1\!-\!\alpha)} \Big) \!+\! \frac{2\epsilon \widetilde{M}}{L(1\!-\!\alpha)}, \nonumber
	\end{align}
	{where} $\widetilde{M} = \sum_{k = 0}^{n - 1} \prod_{s=k+1}^{n-1}  \Big(1 - \frac{\mu_{\sigma(s)}}{L} \Big)$ and $\alpha$ corresponds to the curvature incoherence.
\end{coro}
To elaborate, the permutation map $\sigma$ used by the incremental sampling can be viewed as a particular realization of the random permutation of the random reshuffle scheme. In particular, the convergence error term $\overline{M}$ in \Cref{thm: incohe_sc}  corresponds to the average of the convergence errors over all possible random permutations of the data indexes, whereas the convergence error term $\widetilde{M}$ in \Cref{coro: 2}  is determined by the specific permutation map $\sigma$ used. 
Therefore, depending on the quality of the permutation map, the convergence error of SGD under incremental sampling can be either larger or smaller than that of SGD under random reshuffle.

For SGD with random sampling, we obtain the following convergence rate under minimizer incoherence and a constant step size.
\begin{proposition}[Random sampling]\label{prop: error_bound}
	Let Assumptions \ref{assum: exist} and \ref{assum: f} hold. Apply SGD with random sampling with learning rate $\eta=\frac{1}{L}$ to solve the problem (P). Then, for all $k=nB, n=1,2,...$, 
	\begin{align}
	\mathbb{E} 	\dist_{\Theta^*}^2(\theta_{nB}) \le \big(1 - \frac{\overline{\mu}}{L} \big)^{nB} \Big(	\dist_{\Theta^*}^2(\theta_{0}) - \frac{2 {\epsilon}}{\overline{\mu}} \Big) + \frac{2 {\epsilon}}{\overline{\mu}}. \nonumber
	\end{align}
	where $\bar{\mu} := \frac{1}{n}\sum_{i=1}^n \mu_i$.
\end{proposition}

Comparing the above result with that in \Cref{thm: incohe_sc}, one can see that under minimizer incoherence, SGD with random reshuffle has a better convergence rate coefficient $\alpha=\prod_{i=1}^n (1-\frac{\mu_i}{L})$  than that $(1 - \frac{\overline{\mu}}{L} )^{n} $ of SGD with random sampling (due to the AM-GM inequality). Moreover, regarding the convergence error, one can show that the convergence error $\frac{2\epsilon \overline{M}}{L(1-\alpha)}$ of SGD with random reshuffle is smaller than that $\frac{2{\epsilon}}{\overline{\mu}}$ of SGD with random sampling, and we outline the proof below.
\begin{align}
	\frac{2\epsilon \overline{M}}{L(1-\alpha)} &= \frac{2\epsilon \sum_{k = 0}^{n - 1} \mathbb{E}_{\xi} \big[ \prod_{s=k+1}^{n-1}  \big(1 - \frac{\mu_{\xi(s)}}{L} \big) \big]}{L(1-\alpha)} \nonumber\\
	&\overset{(i)}{\le} \frac{2\epsilon \sum_{k = 0}^{n - 1} \big(1 - \frac{\overline{\mu}}{L} \big)^{n-k-1}}{L(1-\alpha)} \nonumber\\
	&= \frac{2\epsilon L}{\overline{\mu}} \frac{1 -  \big(1 - \frac{\overline{\mu}}{L} \big)^{n}}{L(1-\alpha)} \nonumber\\
	&\overset{(ii)}{\le}  \frac{2\epsilon}{\overline{\mu}}. \nonumber 
\end{align}
To elaborate, consider the quantity $\mathbb{E}_{\xi} \big[ \prod_{s=k+1}^{n-1}  \big(1 - \frac{\mu_{\xi(s)}}{L} \big) \big]$ in $\overline{M}$. Note that for each fixed $k$, the samples $\{\xi(s)\}_{s=k+1}^{n-1}$ are drawn from $\{1,...,n\}$ uniformly at random {\em without replacement} due to the {random reshuffle} scheme, and hence the expectation over $\{\xi(s)\}_{s=k+1}^{n-1}$ consists of $\binom{n}{n-k-1}$ number of different combinations. Therefore, the inequality $(i)$ follows from the Maclaurin's inequality. Moreover, the inequality $(ii)$ follows from the AM-GM inequality. Such a comparison result reveals the statistical advantage of random reshuffle over random sampling: random reshuffle visits all data permutations in expectation and leads to a convergence error in spirit of geometric series (i.e., the $\sum\prod$ term in $\overline{M}$), whereas random sampling samples each data uniformly with replacement and leads to a convergence error in spirit of arithmetic mean (i.e., the $\overline{\mu}$ term). 


\subsection{Empirical Verification} 
We verify our theoretical results obtained in this section via experiments on nonconvex phase retrieval. In specific, consider an underlying complex signal $x_0\in \mathbb{C}^d$ with a set of Gaussian measurement vectors $\{a_i \}_{i=1}^m$. The nonconvex phase retrieval model is written as
$y_i = |\langle a_i, x_0\rangle|+ n,$
where $\{y_i \}_{i=1}^m$ are the phaseless observations and $n$ denotes a Gaussian random noise. To retrieve the signal based on the phaseless observations and the Gaussian measurement vectors, we aim to solve the  following nonconvex problem.
$$\min_{x \in \mathbb{C}^d}  \frac{1}{2n}\sum_{i=1}^n \big( y_i -  |a_i^\intercal x |  \big)^2.$$
Due to noise corruption, the sample losses do not share a minimizer and hence have minimizer incoherence. In particular, the minimizer incoherence increases as the noise level increases. Specifically, we generate $x_0$, $a_i$, and $y_i$ from normal distribution with $d=128, n=512$. We repeat each experiment for $300$ times and use learning rate $\eta=0.1$.

\begin{figure}[bth]\centering
	\includegraphics[width=0.98\linewidth]{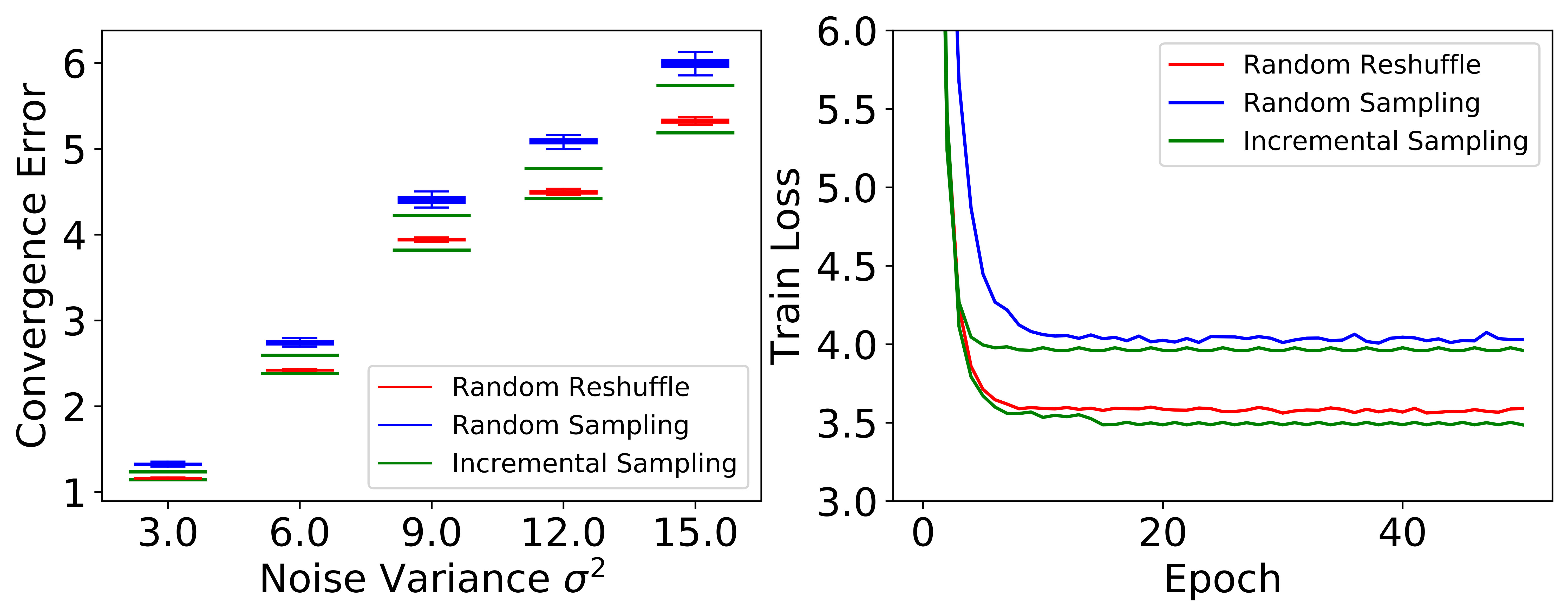} 
	\caption{Left: Impact of minimizer incoherence on convergence error of SGD. Right: Convergence curves of SGD with different sampling schemes under minimizer incoherence.}\label{fig: 5}
\end{figure}

We first explore how the level of noise (i.e., level of minimizer incoherence) in phase retrieval affects the convergence error of SGD with different sampling schemes. \Cref{fig: 5} (Left) presents the box plot of convergence errors of SGD with random sampling and random reshuffle under different levels of Gaussian noise corruptions. For SGD with incremental sampling, we plot the smallest and largest errors achieved in the repeated experiments.
It can be seen that as the noise increases (i.e., minimizer incoherence increases), the convergence errors of these SGDs increase accordingly, which matches our theoretical characterizations of the convergence error. In particular, it can be observed that SGD with random reshuffle consistently has smaller convergence error than SGD with random sampling. Moreover, SGD with incremental sampling can sometimes outperforms SGD with random reshuffle when the permutation map happen to be good.
\Cref{fig: 5} (Right) shows the training loss curves of these algorithms under noise variance $\sigma^2=9$. Under minimizer incoherence, it can be seen that SGD with random reshuffle and SGD with incremental sampling have a comparable convergence speed (i.e., a comparable slope of the training curves), and both of them converge faster than SGD with random sampling. These empirical results validate our convergence rate results of SGD under minimizer incoherence obtained in this section.


 
\section{Conclusion}
In this paper, we propose a model incoherence framework to study the impact of model incoherence on convergence of SGD. When the model has full minimizer coherence, we prove that SGD with random reshuffle converges to a global minimum deterministically and achieves a faster convergence rate than that of SGD with random sampling. When the sample losses have incoherent minimizers, we further show that SGD with random reshuffle has a smaller convergence error than that of SGD with random sampling. Our results reveal the statistical difference between the two random sampling schemes and characterize the impact of model incoherence on the optimization convergence. In the future work, we will further explore the generalization ability of SGD under different sampling schemes and develop a proper analysis framework for it. 

\newpage
\section*{Acknowledgement}
We greatly thank the anonymous reviewers for providing many valuable feedback that help to substantially improve the quality of the paper. 

\bibliography{reference}
\bibliographystyle{icml2019}

\clearpage
\onecolumn
\appendix

\section{Proof of \Cref{prop: vanish_incoherence}}
We first prove the item 1. If $\epsilon_i = 0$ for certain $i$, we have that
\begin{align}
\sup_{\theta \in \Theta^*}   \ell_i(\theta) =  \ell_i^\ast. \nonumber
\end{align}
Since $\ell_{i}^*$ is the global minimum of $\ell_i$, we conclude from the above equality that $\Theta^* \subset \Theta_i^*$. 

Next, we prove item 2. If $\epsilon_i = 0$ for  all $i$, by item 1 we know that $\Theta^* \subset \Theta_i^*$ for all $i$, and hence $\Theta^* \subset \bigcap_{i=1}^n \Theta_i^*$. Now suppose there exists $\theta \in \bigcap_{i=1}^n \Theta_i^* \setminus \Theta^*$. Then, $\theta$ simultaneously minimizes all the sample losses and must be a minimizer of the total loss, i.e., $\theta\in \Theta^*$, contradiction.

\section{Connection between Minimizer Incoherence and other Loss Conditions}\label{subsec: compare_condition}
The notion of minimizer incoherence is related to other loss conditions that have been studied in the existing literature. We outline their connections in this section. 

$\blacktriangleright$ \textbf{Bounded variance}  \shaocong{\cite{ghadimi2013stochastic}}:
In stochastic optimization, it is standard to assume  that the variance of the stochastic gradients is bounded, i.e., for all $\theta\in \mathbb{R}^d$,
\begin{align}
	\mathbb{E}_{\xi}\|\nabla \ell_\xi(\theta) - \nabla f(\theta) \|^2 \le \sigma^2. \label{eq: variance}
\end{align}
In particular, when the total loss $f$ has a unique minimizer $\theta^*$ and all sample losses $\{\ell_i\}_{i=1}^n$ are $1$-gradient dominated\footnote{$\ell$ is called 1-gradient dominated if $\ell(\theta) - \ell^* \le \|\nabla \ell(\theta)\|^2$.}, the stochastic gradient variance at  $\theta^*$ satisfies
\begin{align}
	\mathbb{E}_{\xi}\|\nabla \ell_\xi(\theta^*) - \nabla f(\theta^*) \|^2 
	&\ge \frac{1}{n}\sum_{i=1}^{n} (\ell_i(\theta^*) - \ell_i^*) \nonumber.
\end{align}
in which the right hand side corresponds to the average minimizer incoherence $\frac{1}{n}\sum_{i=1}^{n} \epsilon_i$. Therefore, minimizer incoherence provides an estimate of the stochastic gradient variance at the global minimum, and is weaker than the uniformly-bounded variance condition in \cref{eq: variance}.

$\blacktriangleright$ \textbf{Second moment condition} \cite{bottou2018optimization}: This condition generalizes the previous bounded variance condition as: for some $C\ge 1$ and all $\theta\in \mathbb{R}^d$,
\begin{align}
	\mathbb{E}_{\xi}  \| \nabla \ell_{\xi}(\theta)\|^2 \leq \sigma^2 + C\| \nabla f(\theta) \|^2. \label{eq: second_m}
\end{align} 
In particular, the bounded variance condition corresponds to the second moment condition with $C=1$. In the special case that $\sigma^2=0$ and all the sample losses are convex, the second moment condition implies that $\nabla \ell_i(\theta^*)=0$ for all $i$ and all $\theta^*\in \Theta^*$, i.e.,  every global minimizer of the total loss also minimizes all the sample losses, which further implies full minimizer coherence.

$\blacktriangleright$ \textbf{Interpolation}  \cite{ma2017power}: This condition assumes that the total loss $f$ has a unique minimizer $\theta^\ast$ such that
\begin{align*}
	\ell_i(\theta^\ast ) = \ell_i^* ~\text{for all}~i=1,...,n.
\end{align*}
It can be viewed a special case of the full minimizer coherence, in which the sample losses can share multiple minimizers.

$\blacktriangleright$ \textbf{Growth condition:} In \cite{tseng1998incremental,schmidt2013fast}, the authors considered a  strong growth condition: for some $C \geq 1$ and all $\theta \in \RR^d$,
\begin{align}
	\max_i \| \nabla \ell_i(\theta) \| \leq C \| \nabla f(\theta) \|. \label{eq: strong_grow}
\end{align}  
When all the sample losses are convex, the above condition implies full minimizer coherence. A relaxed version of this condition has been proposed in \cite{vaswani2018fast} as the weak growth condition, which relaxes the $\max_i$ in \cref{eq: strong_grow} to $\mathbb{E}_i$. 

$\blacktriangleright$ \textbf{Expected smoothness} \cite{gower2019}: This condition generalizes the weak growth condition as: 
for some $L> 0$ all $\theta \in \RR^d$,
\begin{align}
	\mathbb{E}_{\xi} \big[\| \nabla \ell_{\xi}(\theta) - \nabla \ell_{\xi}(\theta^*)\|^2 \big] \le L\big(f(\theta) - f^* \big),
\end{align}
where $\theta^*$ is the unique minimizer of $f$. In the case of full minimizer coherence, \cite{gower2019} proved that expected smoothness implies the weak growth condition.


\section{Proof of \Cref{lemma: 1}}
Consider the $k$-th iteration with sample $\xi(k)$. By smoothness of $\ell_{\xi(k)}$, we obtain that
$$\ell_{\xi(k)}(\theta_{k+1} )\leq \ell_{\xi(k)} (\theta_k) + \langle \theta_{k+1}-\theta_k , \nabla \ell_{\xi(k)}(\theta_k) \rangle + \frac{L}{2}\|\theta_{k+1}  - \theta_k \|^2. $$
On the other hand, by restricted convexity of $\ell_{\xi(k)}$, we have: for all $\omega \in \Theta_{\xi(k)}^*$,
$$\ell_{\xi(k)}(\omega) \geq \ell_{\xi(k)}(\theta_k) + \langle \omega - \theta_k, \nabla \ell_{\xi(k)}(\theta_k) \rangle.$$
Combining the above two inequalities yields that
\begin{align*}
\ell_{\xi(k)}(\theta_{k+1} )&\leq   { \ell_{\xi(k)}(\omega)  + \langle  \theta_k -\omega, \nabla \ell_{\xi(k)}(\theta_k)\rangle } + \langle \theta_{k+1}-\theta_k , \nabla \ell_{\xi(k)}(\theta_k) \rangle + \frac{L}{2}\|\theta_{k+1}  - \theta_k \|^2 \\
&= \ell_{\xi(k)}^*  +  \langle \theta_{k+1}-\omega ,  \nabla \ell_{\xi(k)}(\theta_k)   \rangle + \frac{L}{2}\|\theta_{k+1}  - \theta_k \|^2  \\
&= \ell_{\xi(k)}^*  +  \langle \theta_{k+1}-\omega , {-\frac{1}{\eta}(\theta_{k+1} - \theta_k)} \rangle + \frac{L}{2}\|\theta_{k+1}  - \theta_k \|^2 \\ 
&= \ell_{\xi(k)}^*  + { {\frac{1}{\shaocong{2}\eta}\left[ \| \theta_{k} - \omega \|^2 - \| \theta_{k+1} -\omega \|^2 - \| \theta_{k+1} - \theta_{k} \|^2 \right]}}  + \frac{L}{2}\|\theta_{k+1}  - \theta_k \|^2 \\
&= \ell_{\xi(k)}^*  +  \frac{1}{\shaocong{2}\eta}[ \| \theta_{k} - \omega \|^2 - \| \theta_{k+1} -\omega \|^2 ] - \Big(\frac{1}{\shaocong{2}\eta} - \frac{L}{2} \Big)\|\theta_{k+1}  - \theta_k \|^2.
\end{align*} 
Rearranging the above inequality further yields that: for all $\omega \in \Theta_{\xi(k)}^*$,
\begin{align}
\| \theta_{k+1} -\omega \|^2 \leq \| \theta_{k} - \omega \|^2- \shaocong{2}\eta[ \ell_{\xi(k)}(\theta_{k+1} ) - \ell_{\xi(k)}^*  ]  - \big(1 - \shaocong{\eta L} \big)\|\theta_{k+1}  - \theta_k \|^2. \label{eq: lemma1}
\end{align}   
Choose \shaocong{$\eta\le \frac{1}{L}$}, we conclude that for all $\omega \in \Theta_{\xi(k)}^*$,
\begin{align*}
\| \theta_{k+1} -\omega \|^2 \leq \| \theta_{k} - \omega \|^2- \shaocong{2}\eta\big( \ell_{\xi(k)}(\theta_{k+1} ) - \ell_{\xi(k)}^* \big). 
\end{align*}   



\section{Proof of \Cref{lemma: 3}}

{
	
	Note that by Lemma 1, we have that for all $\omega\in \Theta_{\xi(k)}^*$,
	\begin{align*}
	\|\theta_{k+1} - \omega\|^2 &\le \|\theta_{k} -\omega \|^2  - \eta \big( \ell_{\xi(k)}(\theta_{k+1}) -\ell_{\xi(k)}^* \big) \nonumber\\
	&\le \|\theta_{k} -\omega \|^2. \nonumber
	\end{align*} 
	In the case of full minimizer coherence, we have $\Theta^*\subset \Theta_{\xi(k)}^*$. Therefore,  the above result further implies that: for all $k$ and any fixed $\omega \in \Theta^*$,
	\begin{align*}
	\|\theta_{k+1} - \omega\| \le \|\theta_{k} -\omega \| \le \cdots \le \|\theta_{0} -\omega\|<+\infty, 
	\end{align*} 
	where we have used the fact that both $\Theta^*$ and $\theta_{0}$ are bounded. Further notice that $\|\theta_{k+1}\|\le \|\omega\|+\|\theta_{k+1}-\omega\|$, we conclude that the entire trajectory $\{\theta_{k}\}_k$ is bounded. 
}

\section{Proof of \Cref{prop: 2}}
  We first prove item 1. Note that by \Cref{prop: vanish_incoherence} we have $\Theta^* = \bigcap_{i=1}^n \Theta_i^*$. In the proof of Lemma 1 we have shown in \cref{eq: lemma1} that for any $\omega \in \Theta_{\xi(k)}^*$
\begin{align}
\| \theta_{k+1} -\omega \|^2 &\leq \| \theta_{k} - \omega \|^2- \shaocong{2} \eta[ \ell_{\xi(k)}(\theta_{k+1} ) - \ell_{\xi(k)}^*  ]  - \big(1 - \shaocong{\eta L} \big)\|\theta_{k+1}  - \theta_k \|^2. \nonumber
\end{align}   
We can choose any $\omega \in \Theta^*$  and sum the above bound over the $B$-th epoch to obtain that
\begin{align*}
\| \theta_{n(B+1)} -\omega \|^2 \leq \| \theta_{nB} - \omega \|^2- \shaocong{2}\eta \sum_{k=nB}^{n(B+1)-1} \big(\ell_{\xi(k)}(\theta_{k+1} ) - \ell_{\xi(k)}^* \big) - \sum_{k=nB}^{n(B+1)-1}\big(1 - \shaocong{\eta L} \big)\|\theta_{k+1}  - \theta_k \|^2. 
\end{align*}  
Rearranging the above inequality yields that
\begin{align*}
\sum_{k=nB}^{n(B+1)-1} \Big[\big( \ell_{\xi(k)}(\theta_{k+1} ) - \ell_{\xi(k)}^\ast \big) + \big(\frac{1}{\shaocong{2} \eta} - \frac{L}{2}\big)\|\theta_{k+1}  - \theta_k \|^2 \Big] \leq \frac{1}{\shaocong{2} \eta} \big(\| \theta_{nB} -\omega \|^2 -\| \theta_{n(B+1)} -\omega \|^2\big).
\end{align*}  

Further summing the above bound over the epochs $K=0,...,B-1$ yields that
\begin{align}
\sum_{K=0}^{B-1}\sum_{k=nK}^{n(K+1)-1} \Big[\big( \ell_{\xi(k)}(\theta_{k+1} ) - \ell_{\xi(k)}^\ast \big) + \big(\frac{1}{\shaocong{2} \eta} - \frac{L}{2}\big)\|\theta_{k+1}  - \theta_k \|^2 \Big]\leq \frac{1}{\shaocong{2} \eta}\| \theta_{0} -\omega \|^2 .
\end{align}
Note that $ \ell_{\xi_k} (\theta_{k+1}) - \ell_{\xi_k} ^\ast$ is non-negative, and $\big(\frac{1}{\shaocong{2} \eta} - \frac{L}{2}\big)\|\theta_{k+1}  - \theta_k \|^2$ is also non-negative if we choose \shaocong{$\eta \leq \frac{1}{L}$}. Also, the left hand side of the above inequality is  bounded above for all $B$. Therefore, it implies that $ \ell_{\xi_k}(\theta_{k+1}) -  \ell_{\xi_k} ^\ast \overset{k}{\to} 0$, $\|\theta_{k+1}  - \theta_k\|\overset{k}{\to} 0$. In particular, for all subsequences $\{i(T) \}_T, i=1,...,n$, we have $ \ell_{i}(\theta_{i(T)+1}) -  \ell_{i} ^\ast \overset{T}{\to} 0$. Therefore, by continuity of the sample losses, we conclude that all the limit points of $\{\theta_{i(T)+1}\}_T$ belong to the set $\Theta_i^*$ for all $i$. Since $\|\theta_{k+1}  - \theta_k\|\overset{k}{\to} 0$, we conclude that all the limit points $\frakX_i$ of $\{\theta_{i(T)}\}_T$ belong to the set $\Theta_i^*$ for all $i$, and item 1 is proved.

Next, we prove item 2. It suffices to show that $\frakX_i = \frakX_j$ for all $i\ne j$. Consider any $\omega \in \frakX_i$ with a corresponding subsequence $\theta_{i(T_k)} \overset{k}{\to} \omega$. By the random reshuffle sampling, we have $|i(T_k) - j(T_k)| \le n$ for all $i,j,k$. Also, note that $\|\theta_{k+1}  - \theta_k\|\overset{k}{\to} 0$. We obtain that
\begin{align}
\|\theta_{j(T_k)} - \omega \| \le \|\theta_{j(T_k)} - \theta_{i(T_k)} \| + \|\theta_{i(T_k)}  - \omega\| \overset{k}{\to}  0.
\end{align}
Therefore, we showed that every $\omega \in \frakX_i$ is also in any other $\frakX_j$. In summary, $\frakX_i = \frakX_j = \frakX$.
Moreover, since item 1 shows that $\frakX_i \subset \Theta_i^*$, we further obtain that $\frakX \subset \bigcap_{i=1}^n \Theta_i^*$.

\section{Proof of  \Cref{thm: conv}}
	
We prove it by contradiction. Assume there exists $\omega_1, \omega_2 \in \frakX$ such that $\omega_1 \ne \omega_2$. Let $\theta_{q(k)} \to \omega_1$ and $\theta_{p(k)} \to \omega_2$ be two converged subsequences. Without loss of generality, we can always assume that $p(k) > q(k)$ (if not, simply take a subsequence of $\{p(k)\}_k$ such that this property is satisfied). 

Apply the inequality in Lemma 1 with any $\omega \in \frakX\subset \Theta^*$ and note that $p(k)>q(k)$, we obtain that
\begin{align}
\| \theta_{p(k)} - \omega \|  \leq \| \theta_{q(k)} - \omega\|.
\end{align}
In particular, set $\omega = \omega_1$,	the right hand side of the above inequality converges to $0$ because $\omega_1$ is the unique limit point of $\theta_{q(k)}$ by our choice. Therefore, we conclude that $\omega_1$ is also a limit point of $\{ \theta_{p(k)} \}_{k}$, and hence $\omega_1 = \omega_2$, contradiction.

%

\section{Proof of \Cref{thm: rate_coherent-SC}}
Consider the $k$-th iteration with sample $\xi(k)$. By smoothness of $\ell_{\xi(k)}$, we obtain that
\begin{align*}
\ell_{\xi(k)}(\theta_{k+1} )\leq \ell_{\xi(k)} (\theta_k) + \langle \theta_{k+1}-\theta_k , \nabla \ell_{\xi(k)}(\theta_k) \rangle + \frac{L}{2}\|\theta_{k+1}  - \theta_k \|^2.
\end{align*}
On the other hand, by restricted strong convexity of $\ell_{\xi(k)}$, we have: for all $\omega \in \Theta_{\xi(k)}^*$,
\begin{align}
\ell_{\xi(k)}(\omega) \geq \ell_{\xi(k)}(\theta_k) + \langle \omega - \theta_k, \nabla \ell_{\xi(k)}(\theta_k) \rangle + \frac{\mu_{\xi(k)}}{2}\|\theta_{k} - \omega\|^2.
\end{align}
Combining both inequalities above, we obtain that: for all $\omega \in \Theta^*$,
\begin{align}
\ell_{\xi(k)}(\theta_{k+1} ) &\leq  \ell_{\xi(k)}(\omega)  + \langle  \theta_{k+1} -\omega, \nabla \ell_{\xi(k)}(\theta_k) \rangle   + \frac{L}{2}\|\theta_{k+1}  - \theta_k \|^2 - \frac{\mu_{\xi(k)}}{2}\|\theta_{k} - \omega\|^2 \nonumber\\
&=  \ell_{\xi(k)}(\omega)  + \langle  \theta_{k+1} -\omega, \frac{1}{\eta} (\theta_{k} - \theta_{k+1}) \rangle   + \frac{L}{2}\|\theta_{k+1}  - \theta_k \|^2 - \frac{\mu_{\xi(k)}}{2}\|\theta_{k} - \omega\|^2 \nonumber\\
&= \ell_{\xi(k)}(\omega)  + \frac{1}{\shaocong{2}\eta}(\|\theta_{k} -\omega  \|^2 - \| \theta_{k+1} -\omega\|^2)   - \big(\frac{1}{\shaocong{2}\eta}-\frac{L}{2} \big)\|\theta_{k+1}  - \theta_k \|^2  - \frac{\mu_{\xi(k)}}{2}\|\theta_{k} - \omega\|^2.\nonumber 
\end{align} 

Now let \shaocong{$\eta= \frac{1}{L}$}. We further obtain that: for all $\omega \in \Theta^*$,
\begin{align}
\|\theta_{k+1} - \omega\|^2 &\leq \shaocong{\big(1 -  {\mu_{\xi(k)}\eta}  \big)}\|\theta_{k} -\omega \|^2 - \shaocong{2}\eta \Big(   
\ell_{\xi(k)}  (\theta_{k+1}) - \ell_{\xi(k)}^* \Big) \nonumber \\
&\leq \big(1 - \frac{\mu_{\xi(k)}}{L} \big)\|\theta_{k} -\omega \|^2. \label{eq: 3}
\end{align} 
Telescoping the above inequality over the $B$-th epoch and by sampling with random reshuffle, we conclude that: for all $\omega \in \Theta^*$,
\begin{align*}
\| \theta_{n(B+1)} - \omega\|^2 &\leq \prod_{i=1}^n \big( 1- \frac{\mu_{i}}{L} \big)   \|\theta_{nB} - \omega \|^2 \nonumber\\
&= \alpha \|\theta_{nB} - \omega \|^2. \nonumber
\end{align*}
In particular, choose $\omega = \argmin_{u\in \Theta^*} \|\theta_{nB} - u\|$, the above inequality further implies that
\begin{align}
\dist_{\Theta^*}^2(\theta_{n(B+1)}) \le \| \theta_{n(B+1)} - \omega\|^2 \le \alpha \|\theta_{nB} - \omega \|^2 = \alpha	\dist_{\Theta^*}^2(\theta_{nB}). \nonumber
\end{align}
The desired result follows by telescoping the above inequality over the epoch index $B$.

\section{Proof of \Cref{prop: rate_coherent_SGD}}
One can check that \cref{eq: 3} still holds for SGD with random sampling, i.e., 
$$\|\theta_{k+1} - \omega\|^2  \leq \big(1 - \frac{\mu_{\xi(k)}}{L} \big)\|\theta_{k} -\omega \|^2.$$
Taking expectation on both sides of the above inequality yields that
$$\EE \|\theta_{k+1} - \omega\|^2  \leq \big(1 - \frac{ \bar{\mu} }{L} \big)\EE \|\theta_{k} -\omega \|^2,$$
where $\bar{\mu} := \frac{1}{n}\sum_{i=1}^n \mu_i$. 
Telescoping the above inequality over the $B$ epochs yields that
$$\EE \|\theta_{nB} - \omega\|^2  \leq \big(1 - \frac{ \bar{\mu} }{L} \big)^{nB}\EE \|\theta_{0} -\omega \|^2.$$


%

\section{Proof of \Cref{lemma: 2}}
Consider the $k$-th iteration with sample $\xi(k)$. By smoothness of $\ell_{\xi(k)}$, we obtain that
\begin{align*}
\ell_{\xi(k)}(\theta_{k+1} )\leq \ell_{\xi(k)} (\theta_k) + \langle \theta_{k+1}-\theta_k , \nabla \ell_{\xi(k)}(\theta_k) \rangle + \frac{L}{2}\|\theta_{k+1}  - \theta_k \|^2.
\end{align*}
On the other hand, by restricted strong convexity of $\ell_{\xi(k)}$, we have for  $\omega = \proj{\Theta^*}(\theta_k)$,
\begin{align}
\ell_{\xi(k)}(\omega) \geq \ell_{\xi(k)}(\theta_k) + \langle \omega - \theta_k, \nabla \ell_{\xi(k)}(\theta_k) \rangle + \frac{\mu_{\xi(k)}}{2}\|\theta_{k} - \omega\|^2.
\end{align}
Combining both of the above inequalities, we obtain that
\begin{align}
\ell_{\xi(k)}(\theta_{k+1} ) &\leq  \ell_{\xi(k)}(\omega)  + \langle  \theta_{k+1} -\omega, \nabla \ell_{\xi(k)}(\theta_k) \rangle   + \frac{L}{2}\|\theta_{k+1}  - \theta_k \|^2 - \frac{\mu_{\xi(k)}}{2}\|\theta_{k} - \omega\|^2 \nonumber\\
&=  \ell_{\xi(k)}(\omega)  + \langle  \theta_{k+1} -\omega, \frac{1}{\eta} (\theta_{k} - \theta_{k+1}) \rangle   + \frac{L}{2}\|\theta_{k+1}  - \theta_k \|^2 - \frac{\mu_{\xi(k)}}{2}\|\theta_{k} - \omega\|^2 \nonumber\\
&= \ell_{\xi(k)}(\omega)  + \frac{1}{\shaocong{2}\eta}(\|\theta_{k} -\omega  \|^2 - \| \theta_{k+1} -\omega\|^2)   - \big(\frac{1}{\shaocong{2}\eta}-\frac{L}{2} \big)\|\theta_{k+1}  - \theta_k \|^2  - \frac{\mu_{\xi(k)}}{2}\|\theta_{k} - \omega\|^2.\nonumber 
\end{align} 
Choose \shaocong{$\eta= \frac{1}{L}$} and rearrange the above inequality, we obtain that
\begin{align}
\|\theta_{k+1} - \omega\|^2 &\leq \shaocong{\big(1 - {\mu_{\xi(k)}} \eta \big)}\|\theta_{k} -\omega \|^2 - \shaocong{2}\eta \Big(   
\ell_{\xi(k)}  (\theta_{k+1}) - \ell_{\xi(k)} (\omega) \Big)  \nonumber\\
&\leq \big(1 - \frac{\mu_{\xi(k)}}{L} \big)\|\theta_{k} -\omega \|^2 - \shaocong{2}\eta \Big(\ell_{\xi(k)}  (\theta_{k+1}) - \ell_{\xi(k)}^* + \ell_{\xi(k)}^*-\ell_{\xi(k)} (\omega) \Big) \nonumber\\
&\leq \big(1 - \frac{\mu_{\xi(k)}}{L} \big)\|\theta_{k} -\omega \|^2 - \shaocong{2}\eta \Big(\ell_{\xi(k)}  (\theta_{k+1}) - \ell_{\xi(k)}^* + \ell_{\xi(k)}^*-\sup_{\omega\in \Theta^*} \ell_{\xi(k)} (\omega) \Big) \nonumber\\
&\leq \big(1 - \frac{\mu_{\xi(k)}}{L} \big)\|\theta_{k} -\omega \|^2 + \shaocong{2}\eta \epsilon, \label{eq: 4}
\end{align} 
where the last inequality uses the definition of minimizer incoherence, which is bounded by $\epsilon$. 
 Telescoping the above inequality over the iterations of the $B$-th epoch, we obtain that
\begin{align}
\| \theta_{n(B+1)} - \omega\|^2 \leq \prod_{i=1}^n \Big( 1- \frac{\mu_i}{L}\Big)  \|\theta_{nB} - \omega \|^2 +  \shaocong{2}\eta\epsilon \sum_{k = nB}^{n(B+1) - 1} \prod_{s=k+1}^{n(B+1)-1}  \Big(1 - \frac{\mu_{\xi(s)}}{L} \Big)  , \label{eq: tele}
\end{align}
where we define $\prod_{s=n(B+1)}^{n(B+1)-1}  \Big(1 - \frac{\mu_{\xi(s)}}{L} \Big) = 1$ by default. Note that the above inequality is an epochwise contraction with a bounded error term $\eta\epsilon \sum_{k = nB}^{n(B+1) - 1} \prod_{s=k+1}^{n(B+1)-1}  \Big(1 - \frac{\mu_{\xi(s)}}{L} \Big)$, we conclude that $\| \theta_{n(B+1)} - \omega\|^2$ is bounded for all $B$ and hence $\{\theta_k \}_k$ is bounded.

\section{Proof of \Cref{thm: incohe_sc}}
Note that \cref{eq: tele} further implies that
\begin{align}
\dist_{\Theta^*}^2(\theta_{n(B+1)}) &\le \| \theta_{n(B+1)} - \omega\|^2 \nonumber\\
&\le  \prod_{i=1}^n \Big( 1- \frac{\mu_i}{L}\Big)  \|\theta_{nB} - \omega \|^2 + \shaocong{2} \eta   \epsilon\sum_{k = nB}^{n(B+1) - 1} \prod_{s=k+1}^{n(B+1)-1}  \Big(1 - \frac{\mu_{\xi(s)}}{L} \Big) \nonumber\\
&= \prod_{i=1}^n \Big( 1- \frac{\mu_i}{L}\Big) \dist_{\Theta^*}^2(\theta_{nB}) + \shaocong{2} \eta   \epsilon \sum_{k = nB}^{n(B+1) - 1} \prod_{s=k+1}^{n(B+1)-1}  \Big(1 - \frac{\mu_{\xi(s)}}{L} \Big). \label{eq: 1}
\end{align}
Next, denote $\sigma_B$ as the random shuffle permutation performed in epoch $B$ and define the quantity
\begin{align}
M(\sigma_B) := \sum_{k = n(B-1)}^{nB - 1} \prod_{s=k+1}^{nB-1}  \Big(1 - \frac{\mu_{\xi(s)}}{L} \Big)  . \nonumber
\end{align}
It is clear that $M(\sigma_B)$ is a random variable that depends on the permutation $\sigma_B$. We define its expectation as $\mathbb{E}_\sigma M(\sigma_B) := \overline{M}$, which is a fixed constant for every epoch $B$. Then, taking expectation on both sides of \cref{eq: tele} yields that
\begin{align}
\mathbb{E} 	\dist_{\Theta^*}^2(\theta_{n(B+1)}) \leq \alpha \mathbb{E} 	\dist_{\Theta^*}^2(\theta_{nB}) + \shaocong{2} \eta \epsilon\overline{M}. \nonumber
\end{align}
Rearranging the above inequality further yields that
\begin{align}
\mathbb{E} 	\dist_{\Theta^*}^2(\theta_{n(B+1)}) - \frac{\shaocong{2}\eta\epsilon\overline{M}}{1-\alpha} \le \alpha \Big(\mathbb{E} 	\dist_{\Theta^*}^2(\theta_{nB}) - \frac{\shaocong{2}\eta\epsilon\overline{M}}{1-\alpha} \Big), \nonumber
\end{align}
which, after telescoping over $B$, further gives that: for all $B$,
\begin{align}
\mathbb{E} 	\dist_{\Theta^*}^2(\theta_{nB}) \le \alpha^B \Big(	\dist_{\Theta^*}^2(\theta_{0}) - \frac{\shaocong{2}\eta\epsilon\overline{M}}{1-\alpha} \Big) + \frac{\shaocong{2}\eta\epsilon\overline{M}}{1-\alpha}. \nonumber
\end{align}
Lastly, note that we choose \shaocong{$\eta = \frac{1}{L}$}.

\section{Proof of \Cref{coro: 2}}
The proof is similar to that of \Cref{thm: incohe_sc}. The only difference is that the sampling order of the index $\{\sigma(k) \}_k$ is now deterministic. 

One can check that \cref{eq: 1} is valid for SGD with incremental sampling by replacing $\xi(s)$ with $\sigma(s)$, and we have 
\begin{align}
\dist_{\Theta^*}^2(\theta_{n(B+1)}) &\le \prod_{i=1}^n \Big( 1- \frac{\mu_i}{L}\Big) \dist_{\Theta^*}^2(\theta_{nB}) +  \shaocong{2}\eta   \epsilon \sum_{k = nB}^{n(B+1) - 1} \prod_{s=k+1}^{n(B+1)-1}  \Big(1 - \frac{\mu_{\sigma(s)}}{L} \Big) \nonumber\\
&\overset{\text{def}}{:=} \prod_{i=1}^n \Big( 1- \frac{\mu_i}{L}\Big) \dist_{\Theta^*}^2(\theta_{nB}) +  \shaocong{2}\eta   \epsilon \widetilde{M}. \nonumber
\end{align}
Then, the desired result follows from a standard telescoping over $B$ and \shaocong{$\eta = \frac{1}{L}$}.

\section{Proof of \Cref{prop: error_bound}}
One can check that \cref{eq: 4} still holds for SGD with random sampling and step size \shaocong{$\eta=\frac{1}{L}$}. Taking expectations on both sides of the inequality and simplifying yields that
\begin{align}
\mathbb{E} \dist_{\Theta^*}^2(\theta_{k+1}) &\leq \big(1 - \frac{\overline{\mu}}{L} \big) \mathbb{E}\dist_{\Theta^*}^2(\theta_{k}) + \shaocong{2}\eta {\epsilon}, 
\end{align} 
Rearranging and simplifying  the above inequality yields that
\begin{align}
\mathbb{E} 	\dist_{\Theta^*}^2(\theta_{k+1}) - \frac{\shaocong{2}\eta {\epsilon}}{1-(1-\overline{\mu}/L)} \le  \big(1 - \frac{\overline{\mu}}{L} \big) \Big(\mathbb{E} 	\dist_{\Theta^*}^2(\theta_{k}) - \frac{\shaocong{2}\eta {\epsilon}}{1-(1-\overline{\mu}/L)} \Big), \nonumber
\end{align}
which, after telescoping over $k$, further gives that: for all $k=nB$,
\begin{align}
\mathbb{E} 	\dist_{\Theta^*}^2(\theta_{nB}) \le \big(1 - \frac{\overline{\mu}}{L} \big)^{nB} \Big(	\dist_{\Theta^*}^2(\theta_{0}) - \frac{\shaocong{2} \eta {\epsilon}L}{\overline{\mu}} \Big) + \frac{\shaocong{2} \eta {\epsilon}L}{\overline{\mu}}. \nonumber
\end{align}

\end{document}